\documentclass{IEEEtran4PSCC}
\ifCLASSINFOpdf
\usepackage[pdftex]{graphicx}
\else
\usepackage[dvips]{graphicx}
\fi
\usepackage[cmex10]{amsmath}
\usepackage{xcolor}
\usepackage{hyperref}
\usepackage{mathtools}
\newcommand{\defeq}{\vcentcolon=}

\usepackage{amsfonts}
\usepackage{gensymb}
\usepackage{amsthm}
\usepackage{amssymb}
\usepackage[utf8]{inputenc}
\usepackage[T1]{fontenc}
\usepackage[english]{babel}
\usepackage{subcaption}
\newtheorem{theorem}{Theorem}[section]

\theoremstyle{definition}

\theoremstyle{theorem}
\newtheorem*{remark}{Remark}
\usepackage{mathtools}

\usepackage[nocomma]{optidef}
\usepackage[ruled,vlined]{algorithm2e}
\SetKwComment{Comment}{/* }{ */}
\makeatletter
\let\old@ps@headings\ps@headings
\let\old@ps@IEEEtitlepagestyle\ps@IEEEtitlepagestyle
\def\psccfooter#1{%
	\def\ps@headings{%
		\old@ps@headings%
		\def\@oddfoot{\strut\hfill#1\hfill\strut}%
		\def\@evenfoot{\strut\hfill#1\hfill\strut}%
	}%
	\def\ps@IEEEtitlepagestyle{%
		\old@ps@IEEEtitlepagestyle%
		\def\@oddfoot{\strut\hfill#1\hfill\strut}%
		\def\@evenfoot{\strut\hfill#1\hfill\strut}%
	}%
	\ps@headings%
}
\makeatother

\begin{document}
	%
	\title{Towards Balanced Three-phase Charging: Phase Optimization in Adaptive Charging Networks}

	\author{
		\IEEEauthorblockN{Zixin Ye\IEEEauthorrefmark{1}, Tongxin Li\IEEEauthorrefmark{2} and Steven Low\IEEEauthorrefmark{1}\IEEEauthorrefmark{2}}
		\IEEEauthorblockA{\IEEEauthorrefmark{1}Department of Electrical Engineering \\
		\IEEEauthorrefmark{2}Department of Computing and Mathematical Sciences \\
		California Institute of Technology\\
			\{zyye, tongxin, slow\}@caltech.edu}
			}


	\maketitle
	
	\begin{abstract}
		We study the problem of phase optimization for electric-vehicle (EV) charging. 
		We formulate our problem as a non-convex mixed-integer programming problem whose objective is to minimize the charging loss. Despite the hardness of directly solving this non-convex problem, we solve a relaxation of the original problem by proposing the \textit{PXA algorithm} where ``P'', ``X'', and ``A'' stand for three variable matrices in the formed phase optimization problems. We show that under certain conditions, the solution is given by the PXA precisely converges to the global optimum. In addition, using the idea of model predictive control (MPC), we design the {PXA-MPC}, which is an online implementation of the PXA. Compared to other empirical phase balancing strategies, the PXA algorithm significantly improves the charging performance by maximizing energy delivery, minimizing charging price, and assisting future energy planning. The efficacy of our algorithm is demonstrated using data collected from a real-world adaptive EV charging network (ACN).

	\end{abstract}
	

	\thanksto{\noindent Accepted by the 22nd Power Systems Computation Conference (PSCC 2022).}

	\section{Introduction}
	As a long-lasting problem, phase optimization has been studied in recent decades for distribution networks~\cite{TPUVO, OPBDS, TPOPA, feron_2021}.
	A more balanced power supply can improve power delivery efficiency, reduce power line hazards during transmission, and increase decision flexibility for a variety of ancillary services. A balanced three-phase power flow can be produced by dynamically controlling the power supply for load units from each phase line. Such a control mechanism is inherently compatible with an EV charging network. Indeed, phase optimization for EV charging can produce substantial economic benefits, especially in a large-scale and high-density (LSHD), power-delivery adaptive, and information-sensitive network. A three-phase adaptive charging network (ACN)~\cite{ACN2018} includes all those features that strongly motivate the necessity of phase optimization. A real-world ACN often has a three-phase configuration, equipped with software-defined charging control. Our work relies on the Caltech ACN, a multi-story charging facility at Caltech campus~\cite{li2019acndata,lee2021acnsim}.

	\begin{figure}
	\centering
		\includegraphics[width=0.4\textwidth]{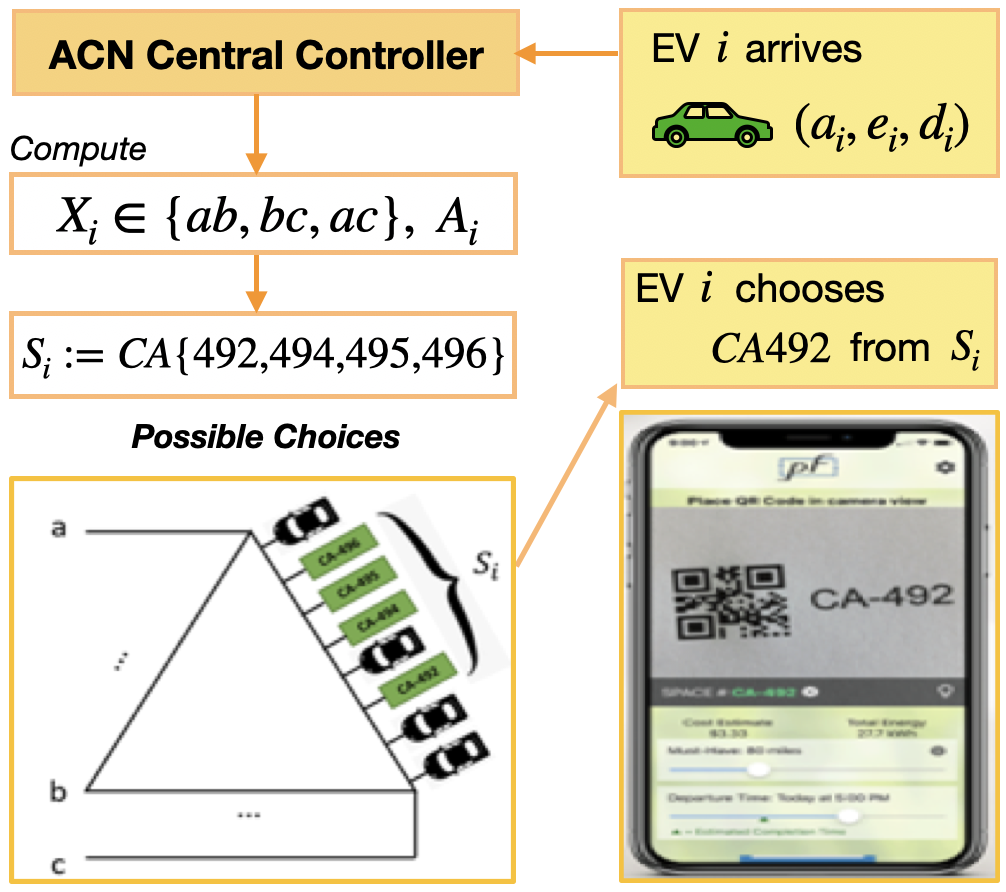}
		\caption{The ACN phase optimization process. \textit{Step 1.} EV $i$ sends information \textit{a priori}: e.g., arrival time $a_i$, parking duration $d_i$, and energy demand $e_i$. \textit{Step 2.} the controller computes the optimal phase selection $X_i\in\{ab,bc,ca\}$ and charging strategy $A_i$. \textit{Step 3.} the controller returns a set $S_i$ of all empty slots whose corresponding phases match $X_i$. \textit{Step 4.} EV $i$ chooses a slot from $S_i$ and starts charging according to $A_i$.}
		\label{fig: PC}
	\end{figure}
	
    In theory, ACN phase optimization is a problem that can be solved by either the demand or the supply side. Existing charging infrastructures, such as~\cite{ACN2018},  assume fixed phase selections determined by demand providers such as EV drivers, and thereby a globally optimal solution can be exactly computed by the supply side. In contrast to this, when the supplier of charging facilities is able to control the phase selections, i.e., the energy demand allocated to phases becomes a decision variable, the charging system can be operated in a more energy and cost-efficient way by maintaining a balanced or near-balanced network. The control of phase selections can be achieved by either deploying PV inverters with a balancing inverter~\cite{LBEVC} or developing smart mobile applications to inform a driver of his desirable parking/charging locations. In fact, balancing phases via either PV inverters or providing phase suggestions to drivers require to solve a phase optimization problem. In this work,  we will formally give a mathematical formulation for the phase optimization problem and provide algorithms for proactive phase allocation.
    
	In practice, we consider an ACN phase optimization process, as illustrated in Figure~\ref{fig: PC}. Based on charging inputs provided by an EV driver upon arrival, an ACN controller computes optimal phase(s) within seconds and shares the suggested location(s) to the driver through a mobile app or payment machines at the entrance of the charging station. Those available locations have chargers mapping to one of the phases, i.e., $\{ab, bc, ca\}$ to which the vehicle will connect. It is worth emphasizing that such a phase-optimization and the demand-allocation process can significantly improve charging efficiency and reduce electricity cost if the following assumptions are satisfied: (1). Each phase is connected with a sufficient amount of ACN chargers. (2). Charger-to-phase connection is immutable. (3) EV drivers' charging inputs are accurate, i.e., the charging duration and energy demand. (4) EV drivers follow suggested charging/parking instructions. Note that (1) and (2) are reasonable assumptions for an LSHD ACN; (3) and (4) are practical assumptions by penalizing the violation of ACN instructions, e.g., by disabling the charging service at places corresponding to the phases that are not allocated to a driver.

	\textit{Related Work.} 
    There is vast literature on the EV charging problem for ACNs. The authors of \cite{li2019acndata} proposed model predictive control-based methods for EV charging. Optimal EV charging as second-order-cone programming (SOCP) was proposed and solved in~\cite{7905971} and \cite{8274166}. The authors of~\cite{lee2020adaptive} illustrated detailed ACN architectures and proposed a novel model-predictive-control (MPC) algorithm that handles unbalanced three-phase infrastructure. With historical ACN charging records, reinforcement learning was used to generate penalty-based terminal functions that facilitate energy-cost efficiency~\cite{li2021learning}.
    Despite that many ACN charging algorithms have been studied and implemented in practice, balancing three-phase EV charging demands for ACNs has not yet attracted enough attention, except for a few works that tackled imbalanced EV connections. For instance,~\cite{EMEV, BCDS} used optimization-based methods to balance renewable penetrations, e.g., to minimize the peak-valley gap. Similarly,~\cite{PVUI, Rezai2016GridablePE} handled imbalanced phase connections by analyzing voltage sensitivity and demand response. However, none of them have studied the phase optimization problem by finding an optimal dispatch of committed EVs for real-world ACNs. 
	
	\textit{Main Contributions.} 
    The main contributions of this paper are three-fold. First, to find optimal phase allocations and charging strategies that minimize a given cost function, we formulate a mathematical model for the ACN phase optimization problem in Section~\ref{pf}. To the best of our knowledge, our work is among the first that considers such an ACN phase optimization problem. Second, we propose an algorithm that finds a near-optimal solution by solving a bi-level optimization and show that it is optimal in the zero laxity case. Moreover, we design an online approach for optimizing and allocating phases in a real-world ACN and validate the efficacy of our algorithms using realistic ACN data.

The rest of this paper is organized as follows. In \autoref{pf}, we describe the model and formulate our phase optimization problem. Next, we demonstrate the intractability of this problem in \autoref{ppop} and propose our method in \autoref{opscad}. Finally, we provide our experimental results in \autoref{sec:exp}.

	\section{Problem Formulation} \label{pf}
	\subsection*{Notation and conventions} 
	\makebox[3cm]{$A_{(i)}, A^{(j)}$}  $i^{th}$ row, $j^{th}$ column of $A$ \par
	\makebox[3cm]{$A \leq B$}  $A_{(i,j)} \leq B_{(i,j)}, \forall i,j$ \par
	\makebox[3cm]{ $\mathbb{R}_{+}, \mathbb{Z}_{+}$}  Non-negative real, integer fields \par
	
	\subsection{EV charging model} \label{evcs}

	Suppose there is a set of $N$ electric vehicles (EVs), denoted by $\mathcal{N} \defeq \{1, \cdots, N\}$ that arrive at an adaptive charging network (ACN) during a discrete time horizon $\mathcal{T} \defeq \{1, \cdots, T\}$.
	Each EV $i$ has a charging-profile $\mathbf{ev}_i \in \mathbb{R}_{+}^3$ defined as a triple $[a_i, d_i, e_i]$ where  $a_i$ (hrs) denotes the EV's arrival time; $d_i$  (hrs) is the charging duration and $e_i$ (kWh) is its initial energy demand. We define the matrix of all $\mathbf{ev}$ as $D \in \mathbb{R}_{+}^{N \times 3}$, where $D_{(i)} \defeq \mathbf{ev}_i$.
Based on $D$, an $N \times T$ binary matrix $E$ with
	\begin{equation}
		\label{eqn:E}
		E_{(i, j)} =
		\begin{cases}
			1      & \quad \text{if } a_i \leq j \leq a_i + d_i\\
			0  & \quad \text{otherwise}
		\end{cases}
	\end{equation} has its
entry $E_{(i,j)}$ indicates whether the EV $i$ is presented at the ACN at time $j$. 
	Moreover, we define an \textit{EV charging matrix} $A \in \mathbb{R}_{+}^{N \times T}$, whose $(i, j)$-th entry is the power (kW) charged for the EV $i$ at time $j$. Finally, we use a $3 \times N$ binary \textit{phase selection matrix} $X$ to denote the assignment of phases for EVs. The set consisting of all such $X$ is 
	\begin{equation}
		\mathcal{X} \defeq \Big\{X \in \{0,1\}^{3 \times N} | \displaystyle \sum ^{3}_{i = 1} X_{(i, j)} = 1, \forall i \in \mathcal{N}\Big\}.
	\end{equation}

In the sequel, with the matrices $D,E,A$ and $X$ and the set $\mathcal{X}$ defined above, we describe the adaptive charging network (ACN) model used in this paper.
	
	\begin{figure}
		\centering
		\includegraphics[width=0.45\textwidth]{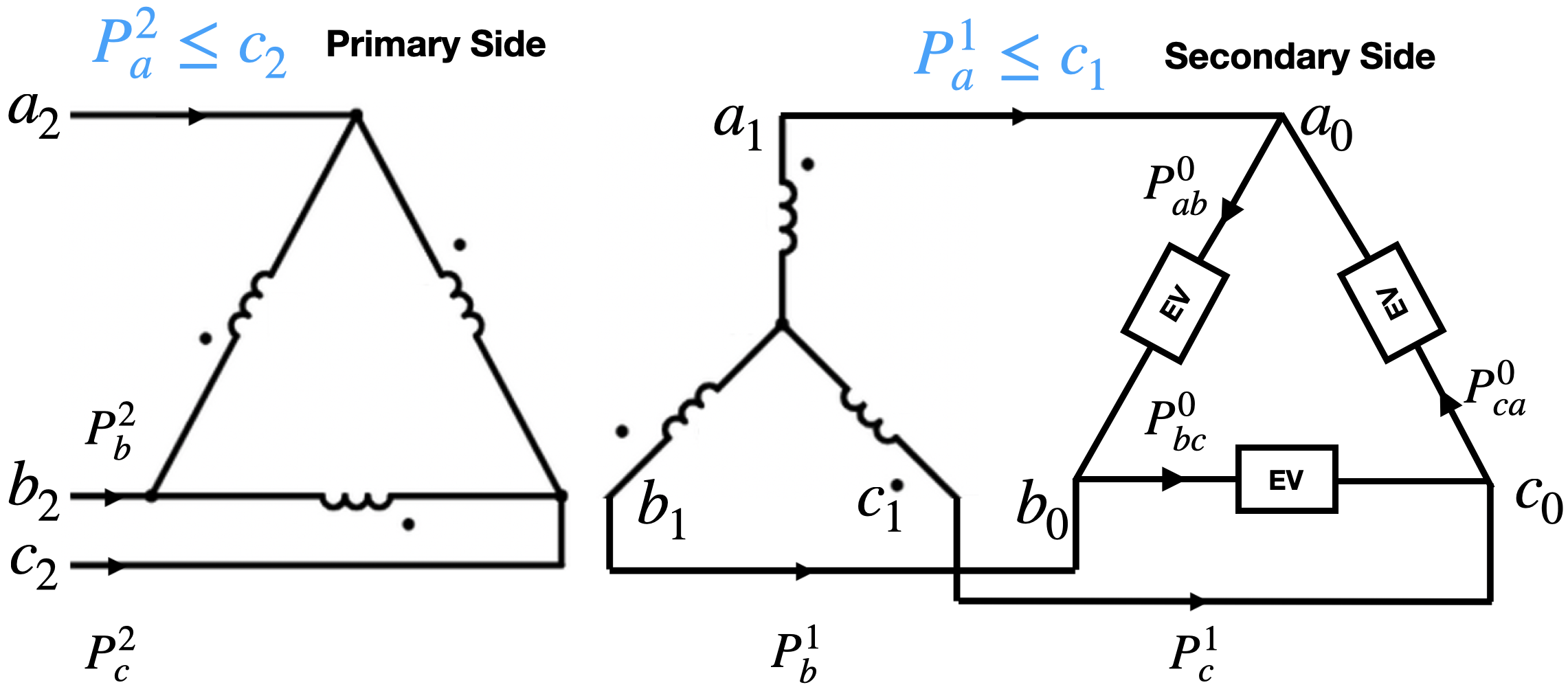}
		\caption{The ACN circuit model involving infrastructure constraints of the Delta-to-Wye step-down transformer and ACN chargers.}
		\label{fig: caltechAcn}
	\end{figure}
	
	\subsection{The ACN model} \label{pfacnf}

	In the following, we describe the constraints of the ACN circuit model depicted in Figure~\ref{fig: caltechAcn}, which has been widely used to model real-world charging facilities such as the California Parking Garage~\cite{7905971,lee2020adaptive}. All charging variables and network constraints are in the unit of power (kW). At each time step $t$, we have the phase power vector  shown in the figure $[P^0_{ab}, P^0_{bc},P^0_{ca}] = (XA)^{(t)}$. Our model considers the charging constraints from the EVSE capacity and EV sessions, along with the network constraint induced by the transformer's capacity. They are the main constraints of the real-world Caltech ACN~\cite{ACN2018}.
	
	\subsubsection{Charging constraints}
First, the charging matrix $A$ must satisfy the charger constraint $\mathcal{C}_{\mathrm{r}}$ such that
	\begin{equation}
		\label{eqn:C_r}
		\mathcal{C}_{\mathrm{r}} \defeq \big\{A \in \mathbb{R}_{+}^{N \times T} | \mathbf{0}_{N \times M} \leq A \leq r_{\max} E\big\}.
	\end{equation}
We assume all chargers have the same maximum power rate $r_{\max} \in \mathbb{R}_{+}$ without loss of generality.
Moreover, $A$ needs to satisfy the demand constraints $\mathcal{C}_{\mathrm{d}}$ such that
	\begin{equation}
		\label{eqn:C_d}
		\mathcal{C}_{\mathrm{d}} \defeq \Big\{A \in \mathbb{R}_{+}^{N \times T} | \Big(\displaystyle\sum_{j = 1}^{T} A^{(j)}\Big)\Delta_T \leq D^{(3)}\Big\}.
	\end{equation}
	
\subsubsection{Network constraints}\label{nc}

The step-down transformer capacity induces line power limits on the primary side $c_2$ and the secondary side $c_1$ in Figure~\ref{fig: caltechAcn}. We denote the set of all $(X,A)\in \mathcal{X} \times \mathbb{R}_{+}^{N \times T}$ satisfying this network constraint by
	\begin{equation}
		\label{eqn:C_soc}
		\begin{aligned}
			\mathcal{C}_{\mathrm{soc}} & \defeq 
			& \left\{\begin{array}{c | c}
				(X, A)  &
				\left|\left[
				\begin{array}{c}
					\Phi_1 \\
					\hline
					\Phi_2
				\end{array}
				\right] X A\right| \leq C_{\max}
			\end{array}\right\}.
		\end{aligned}	
	\end{equation}
In~\eqref{eqn:C_soc}, $|\cdot|$ is the element-wise magnitude operator on the complex matrix, $C_{\max} \in \mathbb{R}^{6 \times T}$ is the line-power limit matrix, and the fraction symbol between $\Phi_1$ and $\Phi_2$ denotes a vertical concatenation of block matrices. The elements of each row depend on the transformer's capacity and e thline power supply. Without loss of generality, mimicking the model used in~\cite{ACN2018,lee2020adaptive}, we consider the effects of the transformer's capacity, i.e. for any $j \in \mathcal{T}$, $C_{\max(i,j)} = c_1$ if $i \in \{1, 2, 3\}$ and $c_2$ if $i \in \{4, 5, 6\}$. The matrix $\Phi_1$ transforms $P_i^0$ to $P_j^1$ and $\Phi_2$ transforms $P_i^0$ to $P_j^2$ ($i \in \{ab, bc, ca\}, j \in \{a, b, c\}$). Given a turning-ratio $n_r$, $\Phi_1$ and $\Phi_2$ can be explicitly defined as
	\begin{equation}
	\label{eqn:Phi1}
	\Phi_1 \defeq 
	\begin{bmatrix}
		1\angle 30 \degree & 0\angle -90 \degree & -1\angle 150 \degree \\
		-1\angle 30 \degree & 1\angle -90 \degree & 0\angle 150 \degree \\
		0\angle 30 \degree & -1\angle -90 \degree & 1\angle 150 \degree
	\end{bmatrix},
    \end{equation}
    \begin{equation}
	\label{eqn:Phi2}
	\Phi_2 \defeq 
	\begin{bmatrix}
		\frac{1}{n_r}\angle 30 \degree & \frac{1}{n_r}\angle -90 \degree & -\frac{2}{n_r}\angle 150 \degree \\
		-\frac{2}{n_r}\angle 30 \degree & \frac{1}{n_r}\angle -90 \degree & \frac{1}{n_r}\angle 150 \degree \\
		\frac{1}{n_r}\angle 30 \degree & -\frac{2}{n_r}\angle -90 \degree & \frac{1}{n_r}\angle 150 \degree
	\end{bmatrix}.
    \end{equation}
	
For notational simplicity, the parameters in our ACN model are summarized as a $4$-tuple: $\mathsf{ACN} \defeq \{D, \mathcal{C}_{\mathrm{r}}, \mathcal{C}_{\mathrm{d}}, \mathcal{C}_{\mathrm{soc}}\}$.
	
	\subsection{Phase optimization for optimal charging}\label{pfocf}
	\label{sec:pop}

	Subject to the constraints described above, the goal of this work is
	to maximize the aggregated energy delivery in $\mathcal{T}$. We define a class of linearly parameterized objective functions $f_{\mathbf{w}}: \mathbb{R}^{N \times T} \to \mathbb{R}$ of the EV charging matrix $A$, where
$
	 	f_{\mathbf{w}}(A) \defeq -\sum^{N}_{i = 1} (A\mathbf{w})_i
$
	 and $\mathbf{w} \in \mathbb{R}_+^T$ is  non-negative.
\begin{remark}
The objective   $f_{\mathbf{w}}$ captures some important properties in the EV charging problem and covers various engineering objectives:
	\begin{enumerate}
		\item if $w_i = 1, \forall i\in \mathcal{T}$, $f_{\mathbf{w}}$ is the power aggregation function. For convenience, we denote such an $f_{\mathbf{w}}$ as $f_{\mathbf{1}}$;
		\item if $w_i = T - i + 1, \forall i \in \mathcal{T}$, $f_{\mathbf{w}}$ is the quick-charge function defined in Eq. 2 of \cite{ACN2018};
		\item if $w_i = p_i \geq 0, \forall i \in \mathcal{T}$, $f_{\mathbf{w}} \Delta_T$ can be used to represent the total energy price paid for the electricity suppliers.
	\end{enumerate}   
\end{remark}

The phase optimization problem considered in this paper is presented below. For the ease of our presentation, we fix $\mathbf{w}=\mathbf{1}$ in the remainder of this paper, and our results can be straightforwardly extended to the general case.
	\begin{mini!}|s|[2]                   
		{\begin{aligned}
				& X \in \mathbb{Z}_{+}^{3 \times N}\\
				& A \in \mathbb{R}_{+}^{N \times T}\\
		\end{aligned}}                               
		{f_{\mathbf{1}}(A) \label{P1:0}}   
		{\label{P:1}}             
		{}                                
		\addConstraint{X}{\in \mathcal{X}, \label{P1:1}}    
		\addConstraint{A}{\in \mathcal{C}_{\mathrm{r}} \cap \mathcal{C}_{\mathrm{d}}, \label{P1:2}}  
		\addConstraint{(X, A)}{\in \mathcal{C}_{\mathrm{soc}}. \label{P1:3}}  
	\end{mini!}
where~\eqref{P1:2} and~\eqref{P1:3} encode the charging and network constraints defined in Section~\ref{pfacnf} respectively.  
For notational convenience,	we define the feasible set as $\mathcal{S}$, the set of all optimal solutions $(X^*, A^*)$ as $\mathcal{S}^*$, the set of all $X^*$ as $\mathcal{X}^*$, and the set containing all $A^*$ as $\mathcal{A}^*$.
	
\section{Preliminaries} \label{ppop}
In this section, we discuss existing approaches for solving the phase optimization formulation~\eqref{P1:0}-\eqref{P1:3}. 
	
\subsection{Existence, non-uniqueness and non-convexity}

	\begin{theorem}[Existence and non-uniqueness]\label{thm:1} The set of optimal solutions $\mathcal{S}^*$ for the optimization~\eqref{P1:0}-\eqref{P1:3} is non-empty and it is a non-singleton set.
	\end{theorem}
	\begin{proof}
		(1) Existence: Given any $X \in \mathcal{X}$, let $\mathcal{A}(X)$ be a set that contains all $A$ such that $A \in \mathcal{C}_{\mathrm{r}} \cap \mathcal{C}_{\mathrm{d}}$ and $(X, A) \in \mathcal{C}_{\mathrm{soc}}$. Since $f_{\mathbf{1}}$ is continuous in $\mathbb{R}^{N \times T}$, it is continuous in $\mathcal{A}(X)$. Moreover, $\mathcal{A}(X)$ is a compact set in $\mathbb{R}^{N \times T}$. According to the Weierstrass extreme value theorem, the infimum of $f_{\mathbf{w}}$ must be attained in $\mathcal{A}(X)$, denoted as $A({X}^*)$. Since this is true for any $X$ and $\mathcal{X}$ is countably finite, $\exists X^* \in \mathcal{X} \text{ s.t. } f_{\mathbf{1}}(A^*({X^*})) \geq f_{\mathbf{1}}(A({X}^*)), \forall X \in \mathcal{X}$. Therefore, $(X^*, A^*({X^*}))$ attains the optimal value. (2) Non-uniqueness: Given any optimal $(X^*, A^*)$, one can verify that any permutation of rows of $X^*\in\mathcal{X}$ is also optimal.
	\end{proof}
    Despite the existence of an optimal solution, finding such a pair $(X^*, A^*)$ is an intractable task because (1) the problem~\eqref{P1:0}-\eqref{P1:3} is an MIP since $X \in \mathbb{Z}_{+}^{3 \times N}$; (2) the set $\mathcal{C}_{\mathrm{soc}}$ is on $(X, A)$. Next, we discuss several existing approaches to solve the optimization~\eqref{P1:0}-\eqref{P1:3} and argue they are not efficient enough for our phase optimization problem.

\subsection{Existing approaches}
\subsubsection{Brute-force interior-point method}\label{socpx}
We first consider a direct approach for solving~\eqref{P1:0}-\eqref{P1:3}. 
	Given a fixed phase selection matrix $X \in \mathcal{X}$, the optimization~\eqref{P1:0}-\eqref{P1:3} is a second-order-cone program (SOCP) formulated in \cite{lee2020adaptive}. Solving the SOCP gives us an exact optimal solution in polynomial time by the interior-point method, according to \cite{Kuo2004}. We denote this SOCP algorithm as $\mathsf{SOCP}({X})$ given $X\in\mathcal{X}$. Since $\mathcal{X}$ is countably finite, a global optimal solution $(X^*, A^*)$ can be exactly found by iterating over all $X \in \mathcal{X}$ in a  brute-force way. We denote this brute-force algorithm as $\mathsf{BFSOCP}$. Due to the brute-force approach,  $\mathsf{BFSOCP}$ is not an efficient method, especially when $N$ and $T$ are large.
	\begin{remark} [EXPTIME]\label{rk:exp}
If the time complexity of $\mathsf{SOCP}$ is $O(p(NT))$, where $p(\cdot)$ is a polynomial function, then the time complexity of $\mathsf{BFSOCP}$ is  $O(3^{p(NT)})$.
	\end{remark}
	
\subsubsection{Approximation by simulated annealing}\label{sec:absa}
Another approach to solving the optimization~\eqref{P1:0}-\eqref{P1:3} is to approximate $(X^*, A^*)$ by simulated annealing, which has been applied to combinatorial problems such as the $q$-coloring problem and the traveling salesman problem. Although this heuristic approach is more efficient than the $\mathsf{BFSOCP}$, it still requires a long time to converge or escape from the local minima. Therefore it becomes unsuitable for the online implementation of our phase optimization problem. In Section~\ref{sec:exp}, we will demonstrate the performance of simulated annealing against other methods.
	
\section{Phase Optimization Algorithms} 
	\label{opscad}

In this section, we present our algorithms for the phase optimization problem defined in \eqref{P1:2}-\eqref{P1:3}. 

\subsection{PXA algorithm for phase optimization}
\label{sofpo}

As discussed in Section~\ref{socpx}, if we have an approach to obtain some optimal phase selection matrix $X^*$ effectively, then we obtain an $(X^*, A^*) \in \mathcal{S}^*$ by solving the $\mathsf{SOCP}({X^*})$ defined in Section~\ref{socpx}. Motivated by this, we propose the following bi-level optimization to generate our estimate of $X^*$, which can be regarded as a charging constraint relaxation of the non-convex optimization~\eqref{P1:0}-\eqref{P1:3}:
	\begin{mini!}|s|[2]                   
		{\begin{aligned}
				& X \in \mathbb{Z}_{+}^{3 \times N}\\
				& P \in \mathbb{R}_{+}^{3 \times T}\\
		\end{aligned}}                               
		{g_{\mathbf{1}}(P) \label{P2:0}}   
		{\label{P:2}}             
		{}                                
		\addConstraint{X}{\in \mathcal{X}, \label{P2:1}}    
		\addConstraint{PM}{\leq XW, \label{P2:2}}  
		\addConstraint{\left|\left[
			\begin{array}{c}
				\Phi_1 \\
				\hline
				\Phi_2
			\end{array}
			\right]P \right|}{\leq C_{\max}. \label{P2:3}}  
	\end{mini!}
In~\eqref{P2:0}-\eqref{P2:3}, the variable $P \in \mathbb{R}_{+}^{3 \times T}$ is  an  \textit{aggregated phase power matrix}, whose row vector is the power charged at $t \in \mathcal{T}$ at each phase. To minimize charging loss as we minimize $f_{\mathbf{1}}$ in~\eqref{P1:0}-\eqref{P1:3}, the aggregated power function $g_{\mathbf{1}}: \mathbb{R}^{3 \times T} \to \mathbb{R}$ is defined as
$
		g_{\mathbf{1}}(P) \defeq -\sum^{N}_{i = 1} \sum^{T}_{j = 1} P_{(i,j)}.
$ Two matrices
	$M \in \mathbb{R}_{+}^{T \times K}$ and $W \in \mathbb{R}_{+}^{N \times K}$ in~\eqref{P2:2} form $K$ linear constraints over $P$ and $X$ as a relaxation of the charging constraint in~\eqref{P1:2}. We will discuss detailed constructions of $M$ and $W$ later in Section~\ref{sec:analysis}. The second-order cone constraint regarding $P$ can be written as~\eqref{P2:3}. We denote the optimal value of the optimization~\eqref{P2:0}-\eqref{P2:3} as $p^*_2$, and the sets containing all feasible $P$ as $\mathcal{P}$ and optimal $P^*$ as $\mathcal{P}^*$.
	\begin{remark}[\textsf{MISOCP}]
The optimization~\eqref{P2:0}-\eqref{P2:3} is a mixed-integer second-order cone programming (MISOCP) since $g_{\mathbf{1}}(P)$ is linear on $P$, ~\eqref{P2:1} is an integral polytope,~\eqref{P2:2} is a mixed-integer polytope on $P$ and $X$ and~\eqref{P2:3} is a second-order cone on $P$. Note that a MISOCP can be efficiently solved by existing optimization solvers. Moreover, we can show that the branch-and-bound method guarantees the convergence to a global optimum with mild conditions.
\end{remark}

By minimizing the variables $P,X$ and $A$ of the bi-level optimizations~\eqref{P2:0}-\eqref{P2:3} and~\eqref{P1:0}-\eqref{P1:3}, we present our \textit{PXA algorithm} in Algorithm~\ref{alg:PXA}. As a MISOCP solver, the computational complexity of PXA is exponential in $N$. However, in practice, the PXA computes much faster than both $\mathsf{BFSOCP}$ and Simulated annealing, especially in the online setting.

\begin{algorithm}[h]
\DontPrintSemicolon
\SetKwBlock{Begin}{function}{end function}
\KwIn{ACN parameters $\mathsf{ACN}=\{D, \mathcal{C}_{\mathrm{r}}, \mathcal{C}_{\mathrm{d}}, \mathcal{C}_{\mathrm{soc}}\}$}
\KwOut{Phase selection $X$ and charging decision $A$}
Initialize $M$ and $W$\\
\Begin($\text{PXA} {(}\mathsf{ACN}|M, W {)}$){
\textit{Step 1}: Solve~\eqref{P2:0}-\eqref{P2:3} with $M, W$ and obtain $X$\\
\textit{Step 2}:
Solve $\mathsf{SOCP}({X})$ and obtain $A$}
\caption{{PXA} for Phase Optimization}
\label{alg:PXA}
\end{algorithm}

	\subsection{MPC-based PXA Algorithm for online phase optimization }\label{sopo}
In real-world EV charging applications, we need to implement the PXA algorithm shown in Algorithm~\ref{alg:PXA} in an online setting, as presented below.

	 Given any $t \in \mathcal{T}$, we observe a fleet of EVs which has already been connected to the ACN (old) and a group of incoming EVs (new). Denote the number of those old and new EVs by $N_{\mathsf{old}}, N_{\mathsf{new}}$. The number of all observable EVs at time $t$ satisfies $N_{\mathsf{old}} + N_{\mathsf{new}} \leq N$. We denote the fixed phase selection matrix for those connected EVs by $X_{\mathsf{old}} \in \mathbb{R}^{3 \times N_{\mathsf{old}}}$. Then $X_{\mathsf{new}} \in \mathbb{R}^{3 \times N_{\mathsf{new}}}$ for those new EVs as our decision variables. 
	 We also define the \textit{EV-profile matrix} at time $t$ as $D(t)$ such that $D(t)^{(1)} = \max(t\mathbf{1}_{N\times 1}, D^{(1)})$, $D(t)^{(2)} = \min(D^{(1)} + D^{(2)} - t, D^{(2)})$ and $D(t)^{(3)} = D^{(3)} - \Delta\mathbf{e}(t)$, where $\Delta\mathbf{e}(t)$ is a vector whose entry $i$ is the energy sum of EV $i$ that has been charged in $\{1,\cdots, t - 1\}$. This online problem is essentially the same as the optimization~\eqref{P2:0}-\eqref{P2:3} except that a subset of columns of $X$, i.e., $X_{\mathsf{old}}$ has been fixed. 
	 Using the idea of MPC, the online MPC-based PXA algorithm ({PXA-MPC}) described above is presented in Algorithm~\ref{alg:1}.
	 \begin{algorithm}
	 	\caption{{{PXA-MPC}} (Online Phase Optimization)}\label{alg:1}
	 	\KwIn{$\mathsf{ACN}= \{D, \mathcal{C}_{\mathrm{r}}, \mathcal{C}_{\mathrm{d}}, \mathcal{C}_{\mathrm{soc}}\}$}
	 	\KwOut{Phase selections $\{X(t)|t\in\mathcal{T}\}$ and charging decisions $\{A(t)|t\in\mathcal{T}\}$}
	 	Initialize 
	 	$A^*(0) \gets \mathbf{0}_{N \times T}$, 
	 	$X(0) \gets \mathbf{0}_{3 \times N}$, $M$ and $W$\\
	 	\For{$t \in \mathcal{T}$}{

	 		\textit{Step 1}: Solve~\eqref{P2:0}-\eqref{P2:3} with  $M$, $W$, $D(t)$, $X_{\mathsf{old}}$, $\mathsf{ACN}$ and obtain  $X_{\mathsf{new}}$
	 		
	 		\textit{Step 2}: Update $X(t) = [X_{\mathsf{old}} |  X_{\mathsf{new}}] $
	 		
	 		\textit{Step 3}: Solve $\mathsf{SOCP}({X})$ with $D(t), X(t), \mathsf{ACN}$. Obtain a solution $A(t)$
	 		
	 		\textit{Update system states}:
	 		
	 		$D^{(2)} \gets \min(D^{(1)} + D^{(2)} - t - 1, D^{(2)})$
	 		
	 		$D^{(1)} \gets \max((t + 1)\mathbf{1}_{N\times 1}, D^{(1)})$
	 		
	 		$D^{(3)} \gets D^{(3)} - {A(t)}^{(1)} \Delta_T$
	 	}
	 \end{algorithm}

	\section{Theoretical Results}
	\label{sec:analysis}
In this section, we provide an analysis of the PXA algorithm, design of the matrices $M$ and $W$, and show that the solution given by the PXA algorithm is optimal under certain laxity conditions.

	\subsection{Constructions of matrices $M$ and $W$}\label{nzlcs}
We propose the following constructions of $M$ and $W$ as a charging constraint relaxation to approximate an optimal $X^*$. 

\begin{figure}
	\centering
		\includegraphics[width=0.43\textwidth]{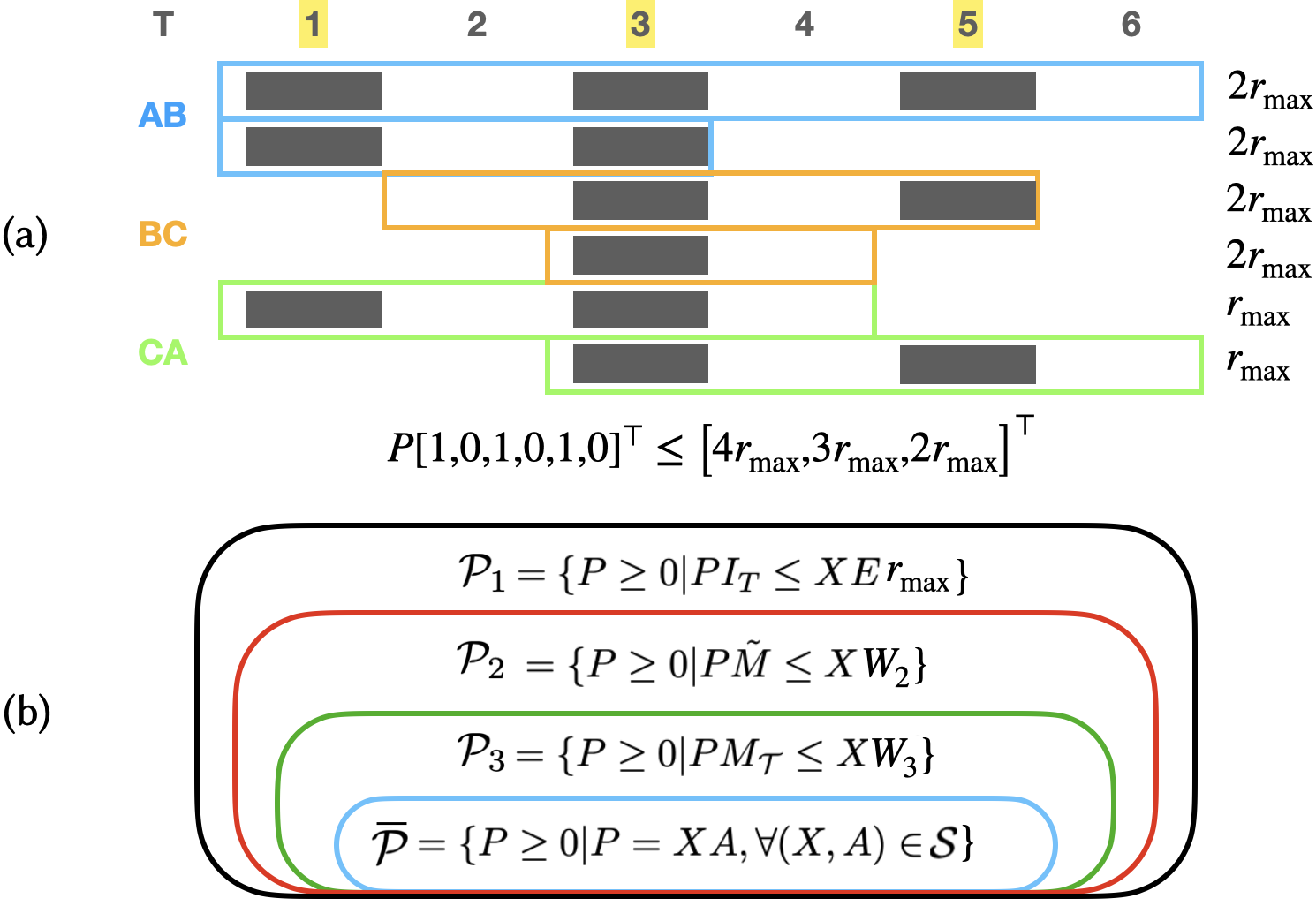}
		\caption{(a) shows an example of charging constraint formed by summing up the phase power at time step $t=1,3$ and $5$. Every two EV sessions are allocated at one phase and their energy demand equals to a proportion of $r_{\max}$. (b) shows the relations of the feasible sets $\mathcal{P}_1,\mathcal{P}_2$ and $\mathcal{P}_3$ formed by different constructions of $M$.}
		\label{fig:pxa}
	\end{figure}

Define $M(m,T) \in \{0, 1\}^{n \times {T \choose m}}$ a  binary \textit{$(m,T)$-selection matrix} if all column vectors of $M(m,T)$ form all possible combinations of selecting $m$ out of $T$ entries. If the $i$-th entry is selected by the $j$-th column, $M_{(i,j)}(m,T)=1$ and $M_{(i,j)}(m,T)=0$ otherwise.
A binary $(m,T)$-selection matrix is used to represent the aggregated phase power for time steps in $\mathcal{T}$. For instance, as shown in Figure~\ref{fig:pxa}(a) with $T = 6$, a charging constraint of $P$ can be described with a column vector $[1, 0, 1, 0, 1, 0]^{\top}$ of $M(3,6)$. Similarly, for all $m\in \mathcal{T}$, we will use every column of $M(m,T)$ to describe a charging constraint in the optimization~\eqref{P2:0}-\eqref{P2:3}. Denote the concatenation of $M(m,T)$ by $M_{\mathcal{T}}$:
	\begin{equation}
\nonumber
		M_{\mathcal{T}} \defeq \left[	M(1,T) | M(2,T) | \cdots | M(T,T) \right] \in \{0, 1\}^{T \times 2^{T}-1}.
	\end{equation}
Moreover, define $2^{T}-1$ horizontal concatenations of power demand vector $D^{(3)}(\Delta_T^{-1})$ as $V \in \mathbb{R}^{N \times 2^{T}-1}$ such that $V\defeq [D^{(3)}(\Delta_T^{-1}) | \cdots | D^{(3)}(\Delta_T^{-1})]$. Denote $C = \min\{A, B\}$ as the minimum-element operation on two matrices $A$ and $B$, e.g., for any $i,j, C_{(i,j)} = \min\{A_{(i,j)}, B_{(i,j)}\}$. Now we can define a set of $P$ as
	\begin{equation}
	\nonumber
		\begin{aligned}
			&\mathcal{P}_3 \defeq \{P \in \mathbb{R}_{+}^{3 \times T} | PM_{\mathcal{T}} \leq XW_3\}, \text{ where }\\
			&   W_3 \defeq \min\{EM_{\mathcal{T}}r_{\max}, V\}
		\end{aligned}
	\end{equation}
so that the aggregated power must be bounded by both the sum of chargers' maximum rates and the total energy demand of EVs that exist in the selected time steps. 
\begin{remark}
We note the following relations:
	\begin{enumerate}
		\item The maximum-rate constraint of $P$ from ACN chargers, i.e. $PI_T \leq XEr_{\max}$, is a superset of $\mathcal{P}_3$, since every column of $I_T$ is a column of $M_{\mathcal{T}}$; 
		\item For every feasible point $(X, A)$ in the optimization~\eqref{P1:0}-\eqref{P1:3}, $(X, P)$ is feasible in $\mathcal{P}_3$, where $P = XA$. $(X, P)$ also satisfies the second-order cone constraint, thus being feasible in the optimization~\eqref{P2:0}-\eqref{P2:3}.
	\end{enumerate}
\end{remark}

Despite that characterizing the constraint~\eqref{P2:2} with $P \in \mathcal{P}_3$ forms an acceptable relaxed constraint as it considers both the charger's maximum rate and the EVs' energy demands, $M_{\mathcal{T}}$ has exponentially many ($2^{T}-1$) columns so that $P\in \mathcal{C}_p$ contains more than $2^{T}-1$ linear constraints. To reduce the number of constraints, we can carefully select $m \ll 2^{T}-1$ columns of $M_{\mathcal{T}}$ to form a relaxation of Constraint~\ref{P2:2}. In our experiments,  we set $\tilde{M} \defeq \left[I | E^T | \mathbf{1}_{T \times 1}\right]$ ($m = T + N + 1$) and its corresponding $W_2$ as a sub-matrix of $W_3$ by selecting the corresponding columns in  $W_3$. Figure~\ref{fig:pxa} (b) shows the relation of the feasible sets with different choices of $M$ and $W$. Although we cannot guarantee  $X \in \mathcal{X}^*$ by using $\tilde M$, we demonstrate in Section~\ref{sec:exp} the result given by the PXA algorithm using $\tilde M$ generates a near-optimal solution. Moreover, in the next section, we prove that in the zero laxity case, considering the relaxation $\mathcal{P}_1$ shown in Figure~\ref{fig:pxa}(b), solving the optimization~\eqref{P2:0}-\eqref{P2:3} by setting $M = I_T$ and $W = Er_{\max}$  gives an optimal $X^*\in \mathcal{X}^*$.
	
	\subsection{Zero laxity analysis} 
	\label{sec:zero-laxity}
Consider the case when all EV sessions have zero laxity, i.e., $\Delta_T r_{\max} d_i = D_{(i)}$, for any  $i \in \mathcal{N}$. Note that the zero-laxity case frequently occurs in specific environments, e.g., fast-charging service centers in roadside rest areas or shopping plazas, wherein most EV drivers leave upon or before their EVs complete charging, thus making EV charging sessions almost zero-laxity. Denote by $I_T$ a $T\times T$ identity matrix. The following theorem holds.

	
		\begin{theorem}[Global optimality]\label{thm:l0o}
Suppose $M = I_T$, $W = Er_{\max}$ and $\Delta_T r_{\max} d_i = D_{(i)}$ for any $i \in \mathcal{N}$. The optimal solution $X$ of the optimization~\eqref{P2:0}-\eqref{P2:3} satisfies that $X \in \mathcal{X}^*$.
		\end{theorem}
		
		\begin{proof}
			It is equivalent to prove that, given any $P^* \in \mathcal{P}^*$, there exists a pair $(X, A) \in \mathcal{S}^*$ such that $P^* = XA$. For any $P^* \in \mathcal{P}^*$, there must exist an $X \in \mathcal{X}$ such that the constraint~\eqref{P2:2} is satisfied. We construct a matrix $A \in \mathbb{R}^{N \times T}$ with the following steps: (1) denote the column sum of $X$ as $\mathbf{n} \defeq [n_1, n_2, n_3]^{\top} = \sum^{N}_{i = 1} X^{(i)}$; (2) for all $i \in \mathcal{T}, j \in \mathcal{N}$,
			\begin{equation}
				\label{eqn:A_c}
				A_{(i, j)} =
				\begin{cases}
					\frac{P^*_{(m, j)}}{n_m}      & \ \text{if } X_{(m, i)} = 1 \text{, where } m \in \{1,2,3\},\\
					0  & \ \text{otherwise}.
				\end{cases}
			\end{equation}
			Notice that $P^* = XA$ and $P^*I \leq XEr_{\max}$. For all $m \in \{1,2,3\}$ and $j \in \mathcal{T}$, we can verify that $P^*_{(m, j)}\leq n_m r_{\max}$, implying that $A_{(i, j)} \leq r_{\max}$. Thus, $A \in \mathcal{C}_{\mathrm{r}}$. Since $\mathcal{C}_{\mathrm{r}} = \mathcal{C}_{\mathrm{r}} \cap \mathcal{C}_{\mathrm{d}}$ and the constraint~\eqref{P2:3} implies $(X, A) \in \mathcal{C}_{\mathrm{soc}}$, we conclude that $(X, A)$ is feasible in the optimization~\eqref{P1:0}-\eqref{P1:3}. We still need to show $p^*_1 = f_{\mathbf{1}}(A)$ to complete the proof. We show this by deriving a contradiction and suppose that there exists some $\tilde A$ that is feasible in the optimization~\eqref{P1:0}-\eqref{P1:3} such that $f_{\mathbf{1}}(\tilde A) < f_{\mathbf{1}}(A)$, thus $p^*_1 \neq f_{\mathbf{1}}(A)$. Since $\tilde A$ is feasible, then there exists an $\tilde X \in \mathcal{X}$ such that $(\tilde X, \tilde A) \in \mathcal{C}_{\mathrm{soc}}$. If we denote $\tilde P \coloneqq \tilde X \tilde A$, it implies $\tilde P$ satisfies the constraint~\eqref{P2:3}. In addition, $\tilde A \in \mathcal{C}_{\mathrm{r}}$ leads to $\tilde X \tilde A \leq \tilde X Er_{\max}$, and hence $\tilde P I_T \leq \tilde X Er_{\max}$. Combining both, $(\tilde X, \tilde P)$ is feasible in the optimization~\eqref{P2:0}-\eqref{P2:3} and $g_{\mathbf{1}}(\tilde P) = f_{\mathbf{1}}(\tilde A) < f_{\mathbf{1}}(A) = g_{\mathbf{1}}(P^*)$. However, $P^*$ is an optimal solution that attains $p^*_2$, implying $g_{\mathbf{1}}(\tilde P) < p^*_2$. Therefore, the existence of $\tilde A$ leads to a contradiction. As a result, $f_{\mathbf{1}}(A) = p^*_1$ implies  that $(X, A) \in \mathcal{S}^*$.
		\end{proof}
	Theorem~\ref{thm:l0o} shows that solving the optimization~\eqref{P2:0}-\eqref{P2:3} gives an exact $X\in\mathcal{X}^*$ in the zero-laxity charging scenarios, thus solving the optimization~\eqref{P1:0}-\eqref{P1:3}. The following theorem shows that solving the optimization~\eqref{P2:0}-\eqref{P2:3} by the branch-and-bound method guarantees the convergence to a global optimum with a mild condition  $p^*_1 \Delta_T = -\sum^N_{i=0}D_{(i, 3)}$ that ensures the feasibility of EVs' energy demands.
	\begin{theorem}\label{thm:eos}
	Suppose $\Delta_T r_{\max} d_i = D_{(i)}$ for any $i \in \mathcal{N}$.	If $p^*_1 \Delta_T = -\sum^N_{i=0}D_{(i,3)}$, solving the optimization~\eqref{P2:0}-\eqref{P2:3} with  $M = I_T$, $W = Er_{\max}$ by branch-and-bound attains $p^*_1$.
	\end{theorem}
	\begin{proof}
		Since Theorem~\ref{thm:l0o} proves that the optimal solution $X$ of the optimization~\eqref{P2:0}-\eqref{P2:3} belongs to $\mathcal{X}^*$, it is enough to show that the branch-and-bound always obtains such an $X$. Since $p^*_1 \Delta_T = -\sum^N_{i=0}D_{(3,i)}$, there exists a $(X^*, A^*)$ feasible in the optimization~\eqref{P1:0}-\eqref{P1:3} such that $f_{\mathbf{1}}(A^*) = p^*_1$. From the proof of Theorem~\ref{thm:l0o}, $(X^*, X^*A^*)$ is feasible in the optimization~\eqref{P2:0}-\eqref{P2:3} and $p^*_2 = g_{\mathbf{1}}(X^*A^*)$. Therefore, we have an integral optimal solution when solving the SOCP relaxation at every branching process. Furthermore, the branch-and-bound process on $X$ is equivalent to a ternary-tree search due to $X \in \{0, 1\}^{3 \times N}$ and $\sum^3_{i = 1} X_{(i, j)} = 1$, for any $j \in \mathcal{N}$. Specifically, at each level $j$ of this ternary-tree, the algorithm computes an optimal solution of the SOCP relaxation with fixes the $j$-th EV's phase selection as $ab, bc, ca$ respectively. We can verify that (i) the branching trajectories of such a ternary-tree are towards all $X \in \mathcal{X}$, thus including the trajectory towards $X^*$; (ii) since $p^*_2$ equals to the total power demand, the bounding criteria is always no less than $p^*_2$, which guarantees that the branching trajectory towards $X^*$ will not be pruned. Combining (i) and (ii), we conclude that the branch-and-bound converges to $p^*_2$, which equals to $p^*_1$.
	\end{proof}

  	\section{Experimental Results} \label{sec:exp}
	\subsection{Experiment setups}
	 We use CVXPY~\cite{diamond2016cvxpy} for modelling the optimization. GUROBI~\cite{gurobi} is used as the branch-and-bound solver for the optimization~\eqref{P2:0}-\eqref{P2:3}, and ECOS~\cite{domahidi2013ecos} for solving the SOCP. We use the ACN-data~\cite{li2019acndata}, which is a dataset collected from a real-world ACN at Caltech. The values of ACN parameters in $\mathsf{ACN}=\{D, \mathcal{C}_{\mathrm{r}}, \mathcal{C}_{\mathrm{d}}, \mathcal{C}_{\mathrm{soc}}\}$ are summarized in Table~\ref{tab:1}. We use  $M = \tilde M$ and $W = W_2$ defined in Section~\ref{nzlcs} in our experiments. We split EV sessions overnight and treat all EV sessions in one day as a charging episode. The charging data $D$ is selected in the pre-COVID period, from date 2018 May. 1st to 2018 Jul. 31st ($91$ episodes in total).

	\subsection{Benchmarks}
	 We consider four benchmarks. (1) {(EV RDM.)} The phases decided by the incoming EV drivers as in the ACN database. (2) {(UNI. RDM.)} We draw each phase connection uniformly at random from $\{ab, bc, ca\}$. (3) {(RRB. RDM.)} We apply the round-robin method based on the arrival time. (4) {(WST.)} The most unbalanced phase selection by allocating all sessions to one phase line. (5) {(M.C.S.A.)} The Markov chain simulated annealing method introduced in Section~\ref{sec:absa} with $10000$ iterations per day. 
	 \begin{figure}
	 	\centering
	 	\begin{subfigure}[b]{0.47\textwidth}
	 		\centering
	 		\includegraphics[width=\textwidth]{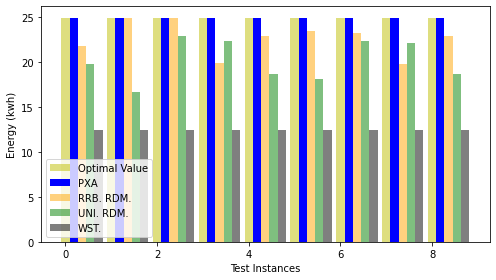}
	 		\caption{}
	 	\end{subfigure}
	 	\hfill
	 	\begin{subfigure}[b]{0.47\textwidth}
	 		\centering
	 		\includegraphics[width=\textwidth]{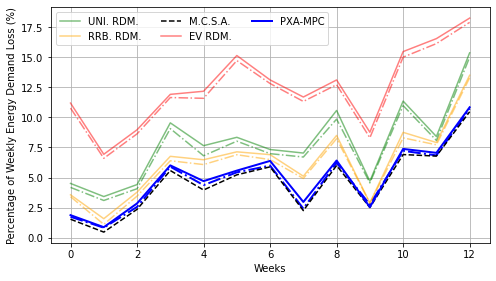}
	 		\caption{}
	 	\end{subfigure}
	 	\hfill
	 	\caption{(a) shows $9$ EV charging instances for various phase optimization algorithms; (b) highlights the performance of the {PXA-MPC} (Algorithm~\ref{alg:1}), compared with the charging outcomes for various phase optimization algorithms.}
	 	\label{fig:onoff}
	 \end{figure}
	 
    \begin{table}[!ht]
    \renewcommand{\arraystretch}{1.3}
    \centering
    \caption{Parameters of the ACN used in our experiments.}
    \label{tab:1}
    \begin{tabular}{|c|c|c|c|c|c|c|}
    \hline
        Episode &  $r_{\max}$& $\mathcal{C}_{\max}$& $V$ & $n_r$& $T$& $\Delta_T$ \\
    \hline
        Day & $3\mathrm{kW}$ & $20\mathrm{kW}$ & $220\mathrm{V}$ & $4$ & $24\mathrm{hrs}$ & $0.2\mathrm{hr}$ \\
    \hline
    \end{tabular}
    \end{table}

\subsection{Near-optimality of the PXA algorithm}
We implement the PXA algorithm in Algorithm~\ref{alg:PXA} for $9$ toy charging scenarios with $T=12$.
The first $4$ scenarios are zero-laxity charging (see Section~\ref{sec:zero-laxity}), and the rest are non-zero laxity charging scenarios. Knowing the ACN parameters $\mathsf{ACN}=\{D, \mathcal{C}_{\mathrm{r}}, \mathcal{C}_{\mathrm{d}}, \mathcal{C}_{\mathrm{soc}}\}$, the PXA algorithm always satisfies $100\%$ of the energy demands, as shown in Figure~\ref{fig:onoff} (a). The uncharged energy demands for {WST.}, {UNI. RDM.}, {RRB. RDM.} are $50\%$, $19.1\%$ and $9.2\%$ respectively.

\subsection{Energy efficiency of the {PXA-MPC} algorithm} 
\begin{figure}
	 	\centering
	 	\begin{subfigure}[b]{0.45\textwidth}
	 		\centering
	 		\includegraphics[width=\textwidth]{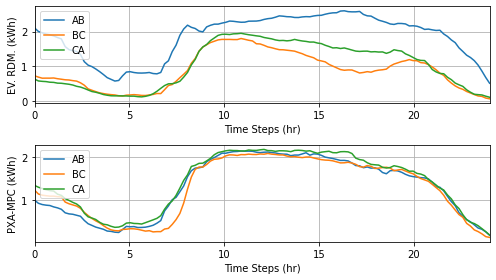}
	 		\caption{}
	 	\end{subfigure}
	 	\hfill
	 	\begin{subfigure}[b]{0.45\textwidth}
	 		\centering
	 		\includegraphics[width=\textwidth]{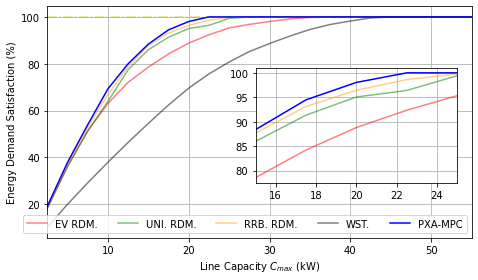}
	 		\caption{}
	 	\end{subfigure}
	 	\hfill
	 	\caption{(a) shows the average energy charged over $91$ episodes at phases $\{ab,bc,ca\}$. The top figure shows the daily trajectory from EV drivers' phase selections and the bottom figure shows the distribution of energy among phases using the {PXA-MPC}; (b) shows demand satisfaction rates with different line-power limit  $C_{\max}$.}
	 	\label{fig:balance}
	 \end{figure}

We implement the {PXA-MPC} (Algorithm~\ref{alg:1}) with the quick-charge function defined in Section~\ref{sec:pop} and assume the controller only gets to know charging sessions when their arrival times are no $\Delta_T$ larger than the current time $t$. In Figure~\ref{fig:onoff} (b), regarding the online performance (solid), the {PXA-MPC} algorithm charges an additional $7.52\%$ of weekly total energy demand compared against the {EV RDM.}; $2.86\%$ compared against the {UNI.RDM} and $1.34\%$ compared against the {RRB. RDM.}. Offline performance (dashed) is also considered. In particular, we observe that the {PXA-MPC} algorithm is near-optimal. Moreover, the {PXA-MPC} takes on average $343.52$ ms to compute at each time step with a standard deviation of $120.42$ ms. In Figure~\ref{fig:balance} (a), we see the {EV RDM.} (top) induces an unbalanced three-phase energy dispatch. While the {PXA-MPC} algorithm (bottom) allocates the energy equally to each phase, especially when the EV demands are high from $10$ am to $8$ pm. We further demonstrate that the {PXA-MPC} algorithm is capable of satisfying more energy demands compared with other phase optimization methods. In Figure~\ref{fig:balance}(b), we show the trade-off between the line capacity and demand satisfaction. We observe that with the same line capacity, the {PXA-MPC} algorithm achieves the highest energy efficiency. In practice, the phases of an ACN can be highly unbalanced, and the {WST.} corresponds to the most unbalanced case when all energy demands are allocated to a single phase.

	\subsection{Electricity cost minimization}
\begin{figure}
		\centering
		\includegraphics[width=0.45\textwidth]{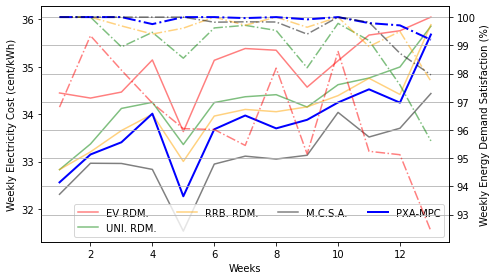}
		\caption{Comparison of  weekly average prices (solid line, left $y$-axis) for {EV RDM.} (red), {UNI. RDM.} (green) and {{PXA-MPC}} (blue), together with their energy delivery in percentage (dashed line, right $y$-axis).}
		\label{fig:price}
	\end{figure}

    Since future information on the energy price is usually unknown, we assume no knowledge of the price vector and use the quick charge function to compute $X^*$ instead. We demonstrate how the {{PXA-MPC}} algorithm can be used to minimize electricity costs. We fix the line capacity as $\mathcal{C}_{\max} = 30$. It can be seen in Figure~\ref{fig:price}, the {{PXA-MPC}} satisfies $99.86\%$ of the energy demand, while the {EV RDM.}, {UNI. RDM.} and the {RRB. RDM.} meet $96.29\%$, $98.85\%$ and $99.51\%$ respectively. Furthermore, the average electricity price of the {{PXA-MPC}} algorithm is $33.85$ cents per kWh, which is constantly lower than that of the {EV RDM.} ($35$ cents/kWh), the {UNI. RDM.} ($34.3$ cents/kWh), and the {RRB. RDM.} ($34.09$ cents/kWh). Compared with the {M.C.S.A}, the {{PXA-MPC}} tends to satisfy higher charging demand in lieu of saving more costs. The phase selection matrix $X^*$ computed from the quick charge function helps the {{PXA-MPC}} satisfy a higher charging demand than that of the {M.C.S.A}. In conclusion, the {{PXA-MPC}} serves more but costs less compared with other online methods.
	\section{Conclusion and Future Directions} \label{c}
		In this paper, we formulate a phase optimization problem for EV charging. We provide a phase optimization algorithm, the PXA algorithm (Algorithm~\ref{alg:PXA}) and an online implementation, the {PXA-MPC} algorithm (Algorithm~\ref{alg:1}). We prove the optimality of the PXA algorithm for the zero-laxity case and demonstrate its efficacy on real-world EV charging data.
Our research opens three interesting directions. Note that our {PXA-MPC} algorithm selects $T + N + 1$ out of $2^T - 1$ linear constraints to solve a linear relaxation of the charging constraint. First, it would be interesting to verify if an optimal solution can be obtained with all linear constraints, i.e., if $\mathcal{P}_3=\overline{\mathcal{P}}$. Second, in addition to the sub-matrix $\tilde M$ defined in Section~\ref{nzlcs}, finding other ways to select those linear constraints and figuring out their optimality would be useful. Finally, it is worth investigating more applications of p

{\footnotesize
    \bibliography{psccBib.bib}{}
    
    \bibliographystyle{ieeetr}
	}
	
\end{document}